\documentclass[11pt]{article}
\usepackage{amsfonts, amsmath, amsthm,amscd,amssymb,verbatim,epsf}
\usepackage{amsmath}
\usepackage{amscd}
\newtheorem{pro}{Proposition}[section]
\newtheorem{thm}[pro]{Theorem}
\newtheorem{lem}[pro]{Lemma}

\newtheorem{cor}[pro]{Corollary}

\theoremstyle{definition}
\newtheorem{dfn}[pro]{Definition}

\theoremstyle{remark}

 \def\d{{\delta}}

 \def\l{{\lambda}}
 \def\L{{\Lambda}}
 \def\m{{\mu}}

 \def\a{{\alpha}}
 
 \def\p{{\partial}}
 
 \def\ra{{\rightarrow}}

 \def\D{{\Delta}}

 \def\z{{\mathbb Z}}
 \def\2{{\mathbb Z_2}}

 \def\sl2{{SL(2,\mathbb C)}}
 
 \def\qed{{\hspace{2mm}{\small $\diamondsuit$}}}

 \def\sl{{{\mbox{\tiny $\L$}}}}

\def\d{\delta}
\def\D{\Delta}
\def\a{\alpha}

\def\L{\Lambda}

\def\l{\lambda}

\begin{document}

\title{Heegaard splittings and virtually Haken Dehn filling II}

\author{Joseph D. Masters, William  Menasco\footnote{{Partially supported by NSF grant \#DMS 0306062}} and
Xingru Zhang\footnote{{Partially supported by NSF grant DMS
0204428.}}}
\maketitle

\begin{abstract}
We use Heegaard splittings to give a criterion for a tunnel number
one knot manifold to be non-fibered and to have large  cyclic
covers. We also show that such a knot manifold (satisfying the
criterion) admits infinitely many virtually Haken Dehn fillings.
 Using a computer, we apply this criterion to
 the 2 generator, non-fibered knot manifolds in the cusped Snappea census.
 For each such manifold $M$, we compute a number $c(M)$,
 such that, for any $n>c(M)$, the $n$-fold cyclic cover of $M$ is large.

\end{abstract}

\section{Introduction}

This paper continues the project, begun in \cite{MMZ},
 of using Heegaard splittings to construct closed essential
 surfaces in finite covers of 3-manifolds.
 The idea, based on the work of Casson and Gordon
 \cite{CG}, is to lift a Heegaard splitting of a 3-manifold
 to a finite cover in which there are disjoint compressing
 disks on each side.  By compressing the lifted Heegaard surface
 along an appropriate choice of such disks,
 we hope to arrive at an essential surface.

 By a \textit{knot manifold} we mean  a connected, compact, orientable 3-manifold whose
 boundary is a single torus.  A \textit{tunnel system}
 for a knot manifold $M$ is a collection $\{t_1, ..., t_n \}$
 where the $t_i$'s are disjoint, properly embedded arcs
 in $M$, such that $M - \mathring{N}(\bigcup t_i)$ is
 homeomorphic to a handlebody.
 The \textit{tunnel number} of $M$, denoted $t(M)$, is the minimal cardinality
 of a tunnel system for $M$.

 We focus attention on tunnel number one,
 non-fibered knot manifolds.  These
 are obtained by attaching a single 2-handle
 to a genus two handlebody.
 We shall give a condition on the 2-handle which,
 if satisfied, ensures that in all large enough cyclic covers,
 the lifted Heegaard splitting can be compressed to obtain
 an essential surface.  There is also a statement
 about incompressibility after Dehn surgery.

 Freedman and Freedman (\cite{FF}) have already proved that,
 in \textit{any} non-fibered knot manifold,
 all but finitely many cyclic covers are large
 (i.e. contain closed essential surfaces).
 Cooper and Long \cite{CL} then proved a result about
 virtually Haken Dehn surgery for these manifolds, and also
 obtained a bound on the number of excluded
 covers in terms of the genus of the knot.

 However, our results provide a computational benefit.
 We computed, for all of the 453
 non-fibered, 2-generator knot manifolds in the SnapPea census,
 a covering degree past which all cyclic covers are large, and
 the bounds obtained are typically improvements over known
 bounds.

  For other connections between Heegaard splittings
 and virtually Haken 3-manifolds, see \cite{L1}, \cite{L2},
 \cite{Li}, \cite{M}.

\section{Definitions, notation, and statement of results}
\label{defs}
 Let $F$ be a connected, closed, orientable surface of positive genus.
 Recall that a \textit{compression body} $W$ is
 a 3-manifold obtained from $(F \times I)$
 by first attaching a collection of 2-handles along disjoint curves
 in one component of $\partial (F \times I)$, and then
 capping off all resulting 2-sphere boundary components with
 3-balls. One of the
 boundary components of $W$ is homeomorphic to $F$,
 and is called the \textit{outer boundary}
 of $W$, denoted $\partial_+ W$.  The other components
 of $\partial W$ form the \textit{inner boundary}, denoted
 $\partial_- W$.

 If $X$ is a 3-manifold with boundary, and $\mathcal{S} \subset X$
 is a collection of disjoint compression disks for $\partial X$, we let
 $X/\mathcal{S} = X - \mathring{N}(\mathcal{S})$.
 If $X$ is a 3-dimensional submanifold of a 3-manifold $Y$, and if
 $\mathcal{S} \subset Y-\mathring{X}$
 is a collection of disjoint compression disks
 for $\partial X$, then we let $X[\mathcal{S}]$ denote
 $X \cup N(\mathcal{S})$, where $N(S)$ is a regular neighborhood of $S$ in
 $Y-\mathring{X}$.

 A \textit{disk system} $\mathcal{S}$ for a compression body $W$ is a set
 of disjoint compressing disks for $\partial W$ of minimal
 cardinality such that $W/\mathcal{S}$ has incompressible boundary.
 We shall use some basic facts about compression bodies,
 which can be found in \cite{Bo}.

 A \textit{Heegaard splitting} of
 a compact 3-manifold is a decomposition
 $M = W_1 \cup_F W_2$, where the $W_i$'s
 are compression bodies with outer boundary
 homeomorphic to $F$.
 The \textit{Heegaard genus} of $M$,
 denoted $g(M)$, is the
 minimal genus of $F$ for all such decompositions.
 If $M$ has boundary, then a \textit{tunnel system}
 for $M$ is a collection of properly embedded arcs
 in $M$, whose exterior is a handlebody.
 The \textit{tunnel number} of
 $M$, denoted $t(M)$, is the minimal cardinality
 among all tunnel systems for $M$.
 It is an elementary fact that, if $M$ is a knot manifold,
 then $g(M) = t(M) + 1$.

 For the remainder of the paper,
 $M$ will be a fixed knot manifold with $t(M) = 1$ and $b_1(M) = 1$.
 Thus there is a Heegaard splitting
 $M = H \cup_F W$, where $H$ is a genus 2 handlebody, and
 $W$ is a genus 2 compression body.
 Let $\mathcal{D} = D_1 \cup D_2$ be a disk system for $H$,
 and let $E$ be the unique (up to isotopy)
 non-separating compression disk for $W$.
 We assume that $E$ has been isotoped so that
 every component of $\partial E - \mathring{N}(\mathcal{D})$ represents an
 essential arc in $F - \mathring{N}(\mathcal{D})$.

 Since $b_1(M) = 1$, there is a unique surjective homomorphism
 $\phi : \pi_1 M \rightarrow \mathbb{Z}$, where $\mathbb Z$ is
 the free factor of $H_1(M)$.
 Let $M_n$ denote the corresponding n-fold cyclic cover,
 with $M_{\infty}$ denoting the infinite cyclic cover.
 Let $H_n, W_n$ and $F_n$ be the pre-images in $M_n$
 of $H, W$ and $F$, respectively.  Then $M_n = H_n \cup_{F_n} W_n$
 is a Heegaard splitting of $M_n$ of genus $n+1$.

 Let $\alpha_1, \alpha_2 \subset F$ be simple closed
 curves transverse
 to $\partial \mathcal{D}$ such that
 $|\alpha_i \cap D_j| = \delta_{ij}$ (the Kronecker delta function).
 We also assume that $\a_1$ and $\a_2$ intersect in a single
 point $p$, which will be the base point for $\pi_1 M$,
 and we assign orientations to $\a_i$ and $\p D_i$
 so that the algebraic intersection numbers $I(\a_i, \p D_i)$ are
 both $+1$. We call such pair of curves $\{\a_1,\a_2\}$ {\it dual curves}
 for the disk system $\{D_1,D_2\}$ of $H$.

\begin{lem}\label{system}
 We may choose a disk system $\mathcal{D}$ so that $\phi( \alpha_1 ) = 0$
 and $\phi( \alpha_2 ) = 1$.
\end{lem}

\begin{proof}
Suppose $\phi( \alpha_i ) = n_i$, and that $|n_1| \leq |n_2|$.
 Let $\d$ be an oriented embedded arc in $\a_1 \cup \a_2$ such that
 $\d \cap (D_1 \cup D_2) = \p \d$, and that its orientation
  agrees with the orientation on $\a_2$,
 but disagrees with the orientation on $\a_1$.
 Let $D_1^{\prime}$ be a properly embedded  disk in $H$ obtained by
   band sum of $D_1$ and $D_2$ along the arc $\d$,
 and let $D_2^{\prime} = D_2$. Then obviously we may assume that $D_1'$ and $D_2'$
 are disjoint, and see that they form  a disk system for $H$.
 Choose  $\a_1^{\prime}=\a_1$, $\a_2^{\prime}=\a_2\a_1^{-1}$. Then up to an
 obvious homotopy   of $\a_2'$ in $\p H$, we may consider $\a_2'$ as a simple closed curve,
 and see that $\{\a_1', \a_2'\}$ form a dual curve pair  for the disk system
 $\mathcal{D}^{\prime}=\{D_1',D_2'\}$ (with a suitable choice of orientation for
 $\p D'_1$). Further
 we have that $\phi(\a_1^{\prime}) = n_1$,
 and $\phi(\a_2^{\prime}) = n_2-n_1$.

 Now we replace the disk system $\mathcal{D}$ with
 the disk system $\mathcal{D}^{\prime}$, and repeat
 the above procedure.
 Applying the Euclidean Algorithm, we may continue until
 we have a disk system $\mathcal{D}$ for which $\phi(\a_1) = 0$ (say)
 and $\phi(\a_2) = gcd(n_1, n_2)$.
 Since $Image(\phi) = \mathbb{Z}$,  $n_1$ and $n_2$ are relatively prime, so the
 resulting disk system satisfies the requirements of Lemma \ref{system}.
 \end{proof}

 For the remainder of the paper, we shall assume
 that the disk system $\mathcal{D}$ has been chosen as
 in Lemma \ref{system}.

 Let $p \in F - \partial \mathcal{D}$ be a base point for $M$,
 and let $p^1, ..., p^n$ be the lifts to $M_n$
 (where $n \in \mathbb{Z}^+ \cup \{ \infty \}$),
 with the natural indexing.
 Let $\delta_1, \delta_2$ in $F - \partial \mathcal{D}$
 be arcs connecting $p$ with $\partial D_i$.
 Let $D_i^j \subset M_n$
 denote the lift of $D_i$ to
 $M_n$ corresponding to the lift of $\delta_i$ with base point $p_j$
  (in our notation, we have suppressed the dependence of $D_i^j$ on $n$,
 trusting the meaning to be clear from context).  See Figure 1.

 Let $I_{geo}(.,.)$ be the geometric intersection pairing, and
 for a loop $\ell$ in $F_{\infty}$,
 define the \textit{width} of $\ell$ to be:
\begin{eqnarray*}
 width(\ell) = Max(j | I_{geo}(\ell, D_2^j) \neq 0)
 - Min(j| I_{geo}(\ell, D_2^j) \neq 0) + 2.
\end{eqnarray*}
 If $\ell$ is a loop in $F$ which lifts to $F_{\infty}$,
 and $\tilde{\ell}$ and $\tilde{\ell}^{\prime}$ are any two lifts
 to $F_{\infty}$, then
 $width(\tilde{\ell})= width(\tilde{\ell}^{\prime})$.
 Thus for any such loop, we define $width(\ell) = width(\tilde{\ell})$,
 where $\tilde{\ell}$ is any lift of $\ell$ to $F_{\infty}$.

Set (for the remainder of the paper)
\begin{eqnarray*}
k = width(\partial E)
\end{eqnarray*}
 We shall restrict attention to the case where
 $E$ intersects $D_1$ and $D_2$ non-trivially.
 If $n \geq k$, let $E^j \subset M_n$ denote the lift
 of $E$ to $M_n$ which intersects $D_1^j$, but is disjoint
 from $D_2^j$.
 Then $D_2^1 \cup \bigcup_{j = 1}^n D_1^j$
 forms a disk system for $H_n$ and $\bigcup_{j=1}^n E^j$
 forms a disk system for $W_n$ (c.f. Fig. 1 for an example
 with $k = 3$ and $n=4$).

\begin{figure}[!ht]
{\epsfxsize=4in \centerline{\epsfbox{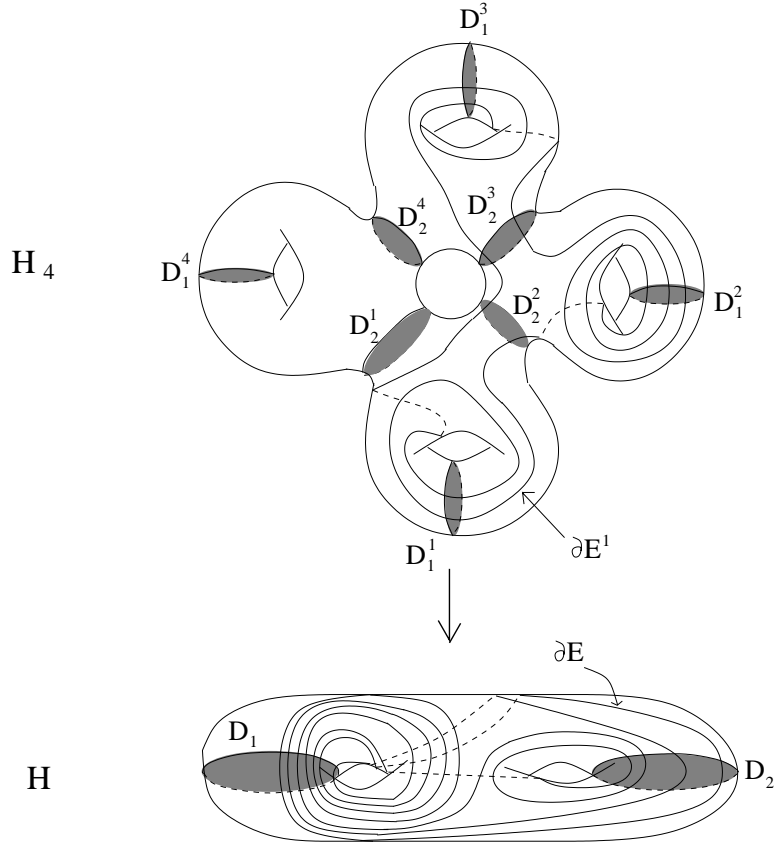}}\hspace{10mm}}
\caption{An example of a 4-fold cyclic cover}
\end{figure}

 Set $H_n^{\prime} = H_n / D_2^1$,
 set  $\mathcal{E}_j = \{ E^1, ...,  E^j \}$,
 and $\mathcal{E}_j^{(i)} = \mathcal{E}_j - E^i$, $1\leq i\leq j$.
Recall that a collection $\mathcal{C}$
 of disjoint simple closed curves in the boundary of a handlebody $X$ is
 \textit{disk busting} if $\partial X - \mathcal{C}$ is
 incompressible in $X$.
\begin{dfn} \label{m-lift}
 Let $m \geq 1$ be an integer.
 We say that $E$ satisfies the \textit{m-lift condition} if:\\
1. The curve $\partial E$ is disk-busting in $H$.\\
2. The set $\partial \mathcal{E}_m$
 is disk-busting in $H_{m+k-1}^{\prime}$.\\
3. The set $\partial \mathcal{E}_m^{(i)}$ is not
 disk-busting in $H_{m+k-1}^{\prime}$, for all $1 \leq i \leq m-1$.\\
4. For each $1 \leq i \leq m-1$, there is a compression disk
 $\Delta_i \subset H_{m+k-1}^{\prime} - \partial \mathcal{E}_m^{(i)}$
 such that $[\partial \Delta_i]$ is linearly independent from
 $\{ [\partial E^{i+1}], ..., [\partial E^m] \}$ in
 $H_1(\partial H_{m+k-1}^{\prime}; \mathbb{Q})$.\end{dfn}
\textit{Remark:} If $M = H[E]$ is fibered, then it follows from
  Lemma \ref{deg2} below that condition 2 fails,
  so $E$ does not satisfy the $m$-lift condition for any $m$.\\

Recall that a 3-manifold $M$ is \textit{large} if it contains a closed
essential surface, i.e. an incompressible surface which is not parallel to a
component of $\partial M$.
We prove:
\begin{thm} \label{cyclic}
 If $E$ satisfies the $m$-lift condition,
 then $M_n$ is irreducible and large for any $n \geq Max(m+k-1, 2k-2)$.
\end{thm}

 Let $\lambda$ be a longitude for $M$ (i.e.
 a simple, closed, essential curve in $\partial M$
 which lifts to a loop in $M_{\infty}$). Fix a meridian $\m$ for $M$ (i.e.
 a simple closed curve in $\p M$ which  intersects
  $\l$ exactly once). A slope $p/q$ in $\p M$ means the pair of homology class
 $\pm(p[\mu]+q[\l])$ and $(p, q)=1$. We use $M(p/q)$ to denote the closed
 manifold obtained by Dehn filling of $M$ with slope $p/q$.
 Let $b=|\phi(\m)|$. Then $b>0$ is a finite integer.    Then we have:

\begin{thm} \label{surgery}
If $E$ satisfies the $m$-lift condition,
   then $M(np/q)$ is virtually Haken for any  $p \geq 2$,
   $n\geq Max\{m+k-1, 2k-2, width(\lambda)+b\}$,
    and q with $(pn,q)=1$.
\end{thm}

Given a 2-generator, 1-relator presentation of a 3-manifold
 group, there is an algorithm to decide if this presentation corresponds
 to a genus 2 Heegaard splitting (conjecturally it always does, see
 \cite{Mas}).
 From the data of such a geometric presentation,
 it is possible to check if the $m$-lift condition holds for a given $m$.
 Using the computer program GAP, we have shown:

\begin{thm}
Every 2-generator, 1-relator 3-manifold $M$ in the SnapPea
 census of 1-cusped hyperbolic 3-manifolds has a genus 2 Heegaard
 splitting. Moreover, if  $b_1(M) = 1$ and $M$ is non-fibered,
 then $M$ has a genus 2 Heegaard splitting
 whose 2-handle satisfies the $m$-lift condition
 for some $m$.
\end{thm}

A complete table of the values of $m$ is available
 at www.math.buffalo.edu/$\tilde{\,\,}$jdmaster.
 The first few values are given in Table 1.

 To prove these theorems, consider the surface $F_n$
 (recall this is the pre-image of the Heegaard surface $F$ in $M_n$).
 Then $F_n$ is a Heegaard surface for $M_n$ of genus $n+1$, which
 we shall compress to both sides.  On the handlebody side, we
 compress $F_n$ along a single lift, $D_2^1$, of $D_2$;
 on the compression body side we compress $F_n$ along all
 the lifts of $E$ which are disjoint from $D_2^1$.
  We shall show that if $E$ satisfies the $m$-lift condition
 and $n \geq m+k-1$, then the resulting surface is incompressible.

Let
\begin{eqnarray*}
 X_n &=& (W_n / \mathcal{E}_{n-k+1}) [D_2^1],\\
 Y_n &=& (H_n/D_2^1)[\mathcal{E}_{n-k+1}],\textrm{ and}\\
 S_n &=& \partial Y_n.
\end{eqnarray*}
 Note that $M_n - \mathring{N}(S_n) \cong X_n \amalg Y_n$.

\begin{lem} \label{nonparallel}
The surface $S_n$ is connected,
 has $genus = k-1$,  and is not parallel to  $\partial M_n$.
\end{lem}

\begin{proof}
To prove that $S_n$ is connected, it is enough to show that
 $[\partial D_2^1], [\partial E^1], ..., [\partial E^n]$  are linearly
 independent in $H_1(F_n)$.

 Let  $I(.,.)$ denote the algebraic intersection pairing on
 the first homology group of a surface.
 Recall from Lemma \ref{system} that $\phi(\a_2)=1$ is a generator
 of $\phi(\pi_1(M))=\z$.  It follows that
 $ I([\partial E], [\partial D_2]) = 0$ and that $\{[\p E], [\p D_2]\}$
 are linearly independent. We may complete
 $\{[\partial E], [\partial D_2] \}$
 to a symplectic basis
 $\{[\partial E], [\partial E]^*, [\partial D_2], [\partial D_2]^* \}$,
 i.e. a basis satisfying
  $I([\partial E], [\partial E]^*) = I([\partial D_2], [\partial D_2]^*) = 1$,
 and
 $I([\partial E]^*, [\partial D_2]) = I([\partial E]^*, [\partial D_2]^*) =
 I([\p E],[\p D_2])=I([\p E],[\p D_2]^*)=0$.
 Let $\alpha \subset F$ be an embedded loop representing
 $[\partial E]^*$, and intersecting
 $\partial E$ geometrically exactly once (such representative always exists).
 Then $\alpha$ lifts homeomorphically to loops
 $\tilde{\alpha}_1, ..., \tilde{\alpha}_n \subset F_n$
 (since $I([\a],[\p D_2])=0$). As $\p E$ also lifts to $F_n$
 as $\p E^1,..., \p E^n$, it is easy to see that
  $I(\tilde{\alpha}_i, \partial E^j) = \delta_{ij}$.
 Also
\begin{eqnarray*}
I([\tilde{\alpha}_i], [\partial D_2^1])
 &=&  I([\tilde{\alpha}_i], [\bigcup_j \partial D_2^j]/n) (\textrm{ since
 all lifts of } \p D_2 \textrm{ are homologous in } F_n) \\
 &=&  I([\alpha], [\partial D_2])/n\\
 &=& 0.
\end{eqnarray*}
 Recall that we have a map $\phi:\pi_1 M \ra \mathbb{Z}$;
 let $\beta \subset \partial M$ be a loop with
 $\phi(\beta) \neq 0$, and let $\beta_n$ be the
 pre-image in $M_n$.
 Then $[\beta_n], [\tilde{\alpha}_1], ..., [\tilde{\alpha}_n]$
 are dual classes for $[\partial D_2^1], [\partial E^1], ..., [\partial E^n]$
 in $H_1(F_n, \mathbb{Q})$,
 which proves the linear independence, and completes the proof that
 $S_n$ is connected.

 The linear independence of
 $[\partial D_2^1], [\partial E^1], ..., [\partial E^n]$
 also allows us to compute:
\begin{eqnarray*}
 genus(S_n) &=& genus(H_n/D_2^1) - |\mathcal{E}_{n-k+1}|\\
 &=& n - (n-k+1) = k-1
\end{eqnarray*}

 Finally, note that $S_n$ is not parallel into $\partial M_n$,
 since every loop in $S_n$ projects to an element in $Ker(\phi)$
 (because $S$ is disjoint from $D_2^1$),
 but each component of $\partial M_n$ contains a loop whose
 projection is not in $ker(\phi)$.
\end{proof}

 To prove Theorems \ref{cyclic} and Theorem \ref{surgery}, we shall show
 that when $n \geq Max\{m+k-1, 2k-2\}$,  $S_n$ is incompressible in both $X_n$
 and $Y_n$, and that when $n \geq Max\{m+k-1, 2k-2, width(\l)\}$, $S_n$
 remains incompressible in an equivariant Dehn filling
 of $M_n$ along $\p M_n$ (which may have several components)
 which is a free cyclic cover of $M(np/q)$.

\section{Background on 1-relator groups and 3-manifolds}

We will require the following result of Magnus (see \cite{MKS}).
The statement given here is
 easily seen to be equivalent to the standard statement.

\begin{thm} \label{frei} (Freiheitsatz for 1-relator groups)
Let
\begin{eqnarray*}
G = <x_1, ..., x_n| w(x_1, ..., x_n)= 1>
\end{eqnarray*}
 be a 1-relator group, where $w$ is a freely reduced word.
 Let $\mathcal{X} = \{ x_1, ..., x_n \}$,
 let $\mathcal{X}^* \subset \mathcal{X}$, and suppose that
 some $x_i \in \mathcal{X}- \mathcal{X}^*$ appears in $w$.
 Then $\mathcal{X}^*$ freely generates a free subgroup of $G$.\qed
\end{thm}

\begin{cor} \label{staggered}
Let $w(x_1, ..., x_k)$ be a word in which
 $x_1$ and $x_k$ appear non-trivially, and consider the group
\begin{eqnarray*}
G= <..., x_{-1}, x_0, x_1, ...
                    |w(x_i, ..., x_{i+k-1})=1, \forall i \in \mathbb{Z}>
\end{eqnarray*}
 Then for any $i$, the set
 $\{ x_i, ..., x_{i+k-2} \}$ freely generates a free subgroup of $G$.
\end{cor}

\begin{proof}
Let $G_i = <x_i, ..., x_{i+k-1}|w(x_i, ..., x_{i+k-1})>$,
 and let $J_i = <x_{i+1}> * \cdots* <x_{i+k-1}>$.
By repeated applications of Theorem \ref{frei},
 the group $G$ has the structure of
 the following iterated amalgamated free product over free subgroups:
 $G \cong \cdots *_{J_{i-1}} G_i *_{J_i} G_{i+1} *_{J_{i+1}} \cdots$.
 By Theorem \ref{frei}, each $J_i$ injects into $G_i$ and $G_{i+1}$.
 Each $G_i$ thus injects into $G$, and so we obtain the corollary.
\end{proof}

Suppose $G = <x,y|w(x,y)>$ is a 1-relator group (where $w$ is a
cyclically reduced word) which admits a surjective homomorphism
 $\psi:G \rightarrow \mathbb{Z}$, such that $\psi(y) = 0$.
 Let $x_i = x^{-i}yx^i$. Then $Ker(\psi)$ is generated by the $x_i$'s,
 and the relation $w$ lifts to a relation $\tilde{w}$ on the $x_i$'s,
 so that $Ker(\psi)$ has a presentation as in the statement of
 Corollary \ref{staggered}.
 Write $\tilde{w} = \Pi_j x_{\mu_j}$, and consider
 the finite integer sequence  $(\mu_j)$. Then we have:

\begin{thm} \label{2gen}(Brown)
 If $(\mu_j)$ has a repeated minimum (or maximum) (i.e. it assumes its minimum
 (or maximum) value more than once), then $ker(\psi)$
 is not finitely generated.
\end{thm}

The case of a repeated minimum is a special instance of Theorem 4.2
 in \cite{B}, and the case of a repeated maximum follows from a trivial
 modification of the proof  (which is an application of the
 Freiheitsatz).

We also need the following, a special case of Corollary 2.2 of \cite{CG}, which
 is in turn a slight modification of a theorem proved by Jaco in \cite{J}.

\begin{thm}(Handle Addition Lemma) \label{handle}
 Let $M$ be an irreducible 3-manifold with compressible
 boundary of genus at least 2, and suppose $\alpha \subset \partial M$ is
 a simple closed curve, such that $\p M - \alpha$
 is incompressible in $M$. Then the 3-manifold obtained by adding a 2-handle
 to $M$ along $\alpha$ is irreducible, and has incompressible boundary.
\end{thm}

\section{2-handles in non-fibered manifolds}
 The results presented in this section are essentially combinations
 of results due to Brown (\cite{B}) and Bieri-Nuemann-Strebel (\cite{BNS}).

If $D$ is a compressing disk for $H_{\infty}$, and $\alpha$ is a
simple  closed curve in $H_{\infty}$, then let $I_{geo}(\alpha, D)$
 be the geometric intersection number of $\alpha$ and $D$;
 in other words, $I_{geo}(\alpha, D)$ is
 the minimal cardinality of $\alpha^{\prime} \cap D$ over
 all curves $\alpha^{\prime}$ which is isotopic  to $\alpha$
 in $H_{\infty}$.

 If $\mathcal{D}$ is a disk system for $H$, then there are
 dual simple loops which form a free basis for $\pi_1 H$.
 If $D \in \mathcal{D}$ corresponds to the generator $x_D$,
 and if $\alpha$ is a simple loop in $H$,
 then $I_{geo}(\alpha, D)$ is the number of times the generator
 $x_D$ appears in a cyclically reduced representative
 for the conjugacy class of $[\alpha]$ in $\pi_1 H$.

\begin{lem} \label{deg2}
 Suppose $M$ is a knot manifold with  $t(M) = 1$, and suppose
 that $E^i$ and $D_i^j$ are
 as defined in Section \ref{defs}.
\newline
(a). If $M$ is fibered, then $I_{geo}(\partial E^1, D_1^1)
                              = I_{geo}(\partial E^1, D_1^k) =1$
in $H_{\infty}$.\newline (b). If $M$ is non-fibered, then
$I_{geo}(\partial E^1, D_1^1) \geq 2$
 and $I_{geo}(\partial E^1, D_1^k) \geq 2$ in $H_{\infty}$.
\end{lem}

\begin{proof}
(a).
 Suppose $I_{geo}(\partial E^1, D_1^1) \geq 2$ or
 $I_{geo}(\partial E^1, D_1^k) \geq 2$.

 Then by Theorem \ref{2gen} (together with the note in
 the proceeding paragraph of
 the present lemma), $Ker(\phi)$ is not finitely generated, and so
 by \cite{S1}, $M$ is not fibered.\\

\noindent (b). Suppose $E^1$ intersects one of the disks,
 say $D_1^k$, exactly once.  We shall show that in this
 case $M$ is fibered.

 Dual to each $D_1^i$ is a generator $x_i$ for the fundamental
 group of $M_{\infty}$.
 The boundary of the disk $E^1$ gives a relation among these generators
 which involves $x_k$ only once; therefore
 $x_k \in <x_1, ..., x_{k-1}> \subset \pi_1 M_{\infty}$.
 Similarly, using the relation corresponding to the disk
 $E^2$, we get that

\begin{eqnarray*}
 x_{k+1} \in <x_1, ..., x_k> = <x_1, ..., x_{k-1}> \subset \pi_1 M_{\infty}.
\end{eqnarray*}

 Continuing in this way, we see that all of the generators $x_i$
 with $i \geq k$ can be expressed in terms of $x_1, ..., x_{k-1}$.

 Let $H^*$ be the component of $H_{\infty}/D_2^1$ containing
 $D_1^1, D_1^2,  ...$, and let $Q = H^*[E^1, E^2, ...]$,
 which is a submanifold of $M_{\infty}$.
 The argument we just gave shows that $\pi_1(Q)$ is finitely generated.

 Note that is a non-separating
  incompressible  surface $S$ in $M$ with boundary slope $\l$
 such that  $M_\infty$ is the infinite cover dual to $S$.
 Let $\tilde{S}$ be a lift of $S$ to $M_{\infty}$
 which is disjoint from $Q$, and let $Q^+$ be the component of
 $M_{\infty} - \mathring{N}(\tilde{S})$ which contains $Q$.
 Then $Q^+ - \mathring{N}(Q)$ is compact, and
 since $\pi_1 Q$ is finitely generated,
 $\pi_1 Q^+$ is finitely generated as well.

 Let $M_0^- = M - \mathring{N}(S)$, let
 $S_0$ and $S_1$ be the two pre-images of $S$ in
 $\partial M_0^-$, let $\tilde{S}_i$ be the pre-images of $S$ in
 $M_{\infty}$, and let $M_i^-$ be the submanifold of $Q$
 bounded by $\tilde{S}_0$ and $\tilde{S}_i$.
 Since $\tilde{S}_i$ is incompressible, $M_i^-$
 is $\pi_1$-injective in $Q$ for each $i$.
 If neither of the maps
 $i_* \pi_1 S_j \rightarrow \pi_1 M_0$
 is onto, then $\{ \pi_1 M_i \}$ forms an
 an infinite sequence of subgroups of $\pi_1 Q$, with $\pi_1 M_{i+1}$
 properly containing $\pi_1 M_i$ for each $i$,
 which is a contradiction, since $\pi_1 Q$ is finitely generated.
 Therefore, one of the induced maps
 $\pi_1 S_j \rightarrow \pi_1 M_0$
 is onto, and so, as in the proof of Theorem 2 in \cite{S1},
 $M$ is fibered.
\end{proof}

\begin{cor} \label{fib}
If $M$ is fibered, then $E$ does not satisfy
 the $m$-lift condition for any $m$.
\end{cor}

\begin{proof}
 Suppose $M$ is fibered, let $m \geq 1$ be an integer, and consider
 the cover $F_{m+k-1}$ of $F$.
 In $F_{m+k-1}$, we have that $E^j$ is disjoint from $D_1^1$
 for all $2 \leq j \leq m-k+1$,
 and by Lemma \ref{deg2} a, $|\partial E^1 \cap \partial D_1^1| = 1$.
  Therefore there is a compressing disk $\Delta$ in $H_{m+k-1}$
 (whose boundary is equal to $\partial N(\partial E^1 \cup \partial D_1^1)$)
 with $\partial \Delta \cap \partial \mathcal{E}_m = \emptyset$, and
  so $E$ fails Condition 2.
\end{proof}

\section{Proof of irreducibility of $M_n$}
We shall now begin the proof of Theorem \ref{cyclic},
 which will occupy the next three sections.
\begin{lem} \label{irred}
If $H - \partial E$ has incompressible boundary, then $M_n$
 is irreducible for all $n$.
\end{lem}

\begin{proof}
 By Theorem \ref{handle},  $M_1 = M$ is irreducible.
 By \cite{MSY} (or \cite{D}), the cover of an
 irreducible manifold is irreducible, so $M_i$
 is irreducible for all $i \geq 2$.
\end{proof}

\section{Proof that $X_n$ has incompressible boundary}

 We remark that
Conditions 2-4 of the $m$-lift condition are not needed in this case,
 and so our proof implies that
 $X_n$ has incompressible boundary whenever $\partial M$ is incompressible.

Let $W_n^{\prime} = W_n / \mathcal{E}_{n-k+1}$.
 By Theorem \ref{handle} it is enough to prove:

\begin{lem} \label{bdy}
If $n \geq 2k-2$, then $W_n^{\prime} - \partial D_2^1$ has incompressible boundary.
\end{lem}

First we need:
\begin{lem} \label{mer}
For each $n$ there is a loop $\a_n \subset \p M_n$ such that
$I(\a_n,D_2^1) \neq 0$.
\end{lem}

\begin{proof}
By the exact sequence of the pair, there is a loop
 $\a \in \p M$ such that $\phi[\a] = I(\a, D_2^1) \neq 0$.
 Letting $\a_n$ be the pre-image of $\a$ in $\p M_n$,
 we have $I(\a_n, \bigcup_{i=1}^n D_2^i) = n I(\a,D_2) \neq 0$.
 Since $[D_2^i] = [D_2^j] \in H_1(M_n)$ for all $i,j$,
 then we have $I(\a_n, D_2^1) = I(\a, D_2) \neq 0$.
\end{proof}

\begin{proof} (of Lemma \ref{bdy})
Suppose otherwise that there is a compressing disk in $W_n^{\prime}
- \partial D_2^1$.
 First, if there is a compressing disk, we claim that there must
 be a non-separating one.
 To see this, suppose that $\D$ is a separating compressing disk.
 If there are no non-separating compressing disks,
 then one of the components of $W_n^{\prime} / \D$
 is homeomorphic to a surface cross an interval, and
 the other component is a handlebody containing $\p D_2^1$.
 Every curve in $\p M_n$ lies on the surface cross interval side,
 but by Lemma \ref{mer}, there is a curve in $\p M_n$ which has non-trivial
 intersection with $D_2^1$, yielding  a contradiction.
 So we may assume that  there is a non-separating compressing disk
 $\D$ in $W_n^{\prime}
- \partial D_2^1$.

 Consider the Heegaard surface $F_n$ and the curves
 $\p D_2^j$ and $\p E^j$, $j=1,...,n$,  in $F_n$.
Note that
  since $n \geq 2k-2$,
 $\{\p E^{n-k+2}, ..., \p E^n\}$ are all disjoint
 from  $\p D_2^{n-k+2}$ and $\p D_2^k$ (by considering
 the definition of  $k=width(E)$).
 The two simple closed curves
$\p D_2^{n-k+2}$ and $\p D_2^k$ cut $F_n$ into two components,
$F_n^1$ and $F_n^2$, one of which, say $F_n^1$, is disjoint from all
$\p E^{n-k+2},...,\p E^{n}$.

The curves $\p E^1,...,\p E^{n-k+1}$ may intersect $\p D_2^{n-k+2}$
and $\p D_2^k$. But by the property of the  width $k$ again, any arc
component of $F_n^2\cap(\cup_{j=1}^{n-k+1} \p E^j)$ has either both
endpoints in $\p D_2^{n-k+2}$ or both endpoints in $\p D_2^k$. Further,
every arc component of $F_n^2\cap(\cup_{j=1}^{n-k+1} \p E^j)$ is
disjoint from $\p D_2^1$. Let $F_n^3$ be the subsurface of $F_n$
which is the union of $F_n^1$ and a small regular neighborhood of
the arcs $F_n^2\cap(\cup_{j=1}^{n-k+1} \p E^j)$ in $F_n^2$ (in other
words, $F_n^3$ is $F_n^1$ with some bands attached, one for each arc
component in $F_n^2\cap(\cup_{j=1}^{n-k+1} \p E^j)$). Then $F_n^3$
is a connected subsurface of $F_n$  which contains all $\p E^1,
...,\p E^{n-k+1}$ but is disjoint from all $\p E^{n-k+2},....,\p
E^n$ and $\p D_2^1$.

Let $F_n^4$ be the surface obtained from $F_n^3$ by surgery along
the curves $\p E^1,...,\p E^{n-k+1}$ (i.e. cut $F_n^3$ open along
$\{\p E^1,...,\p E^{n-k+1}\}$ and fill each of new
 boundary circles
with a disk), which may not be connected. Note that $\{E^{n-k+2},
..., E^n\}$ is a disk system for the compression body $W_n'$ and
$W_n'/\{E^{n-k+2},...,E^n\}$ is an $I$-bundle over a surface. As
$F_n^4$ is disjoint from $\{\p E^{n-k+2}, ..., \p E^n\}$,  $F_n^4$
is contained in the horizontal boundary of the $I$-bundle. So $\p
F_n^4\times I$ are vertical annuli of this $I$-bundle. By standard
cut-and-paste operations along arcs and circles of $\D\cap (\p
F_n^4\times I)$, we get a non-separating compressing disk, still
denoted $\D$, which is disjoint from the annuli $\p F^4_n\times I$.
It is easy to see that  $\D$ cannot be contained in $F_n^4\times I$,
 so it
follows that $\p \D$ is disjoint from all $\p D_2^1$, $\p
D_2^k$,..., $\p D_2^{n-k+2}$. Hence the width of $\D$ is strictly
less than $k$.

 Let $\d = \p \D$, and
 let $\tilde{\d}$ be a lift of  $\d$ in $M_{\infty}$,
The group $\pi_1 M_{\infty}$ has the following presentation:
$$ \pi_1 M_{\infty} = < ..., x_{-1}, x_0, x_1, ... | w(x_i, ..., x_{i+k-1})
 = 1, \, \forall i \in \mathbb{Z}>,$$
where the relations correspond to the lifts of $E$.
 Since $width(\d) < k$, the loop $\tilde{\d}$ represents an element
 in the subgroup of $\pi_1 M_{\infty}$
 generated by the elements $x_i, ..., x_{i+k-2}$, for some $i$.
 By Corollary \ref{staggered}, these elements
 are a basis for a free subgroup;
 since $\tilde{\d}$ represents a trivial element in $\pi_1 M_{\infty}$,
 we see that $\tilde{\d}$ represents the trivial word
 in $x_i, ..., x_{i+k-2}$.
 Thus $\tilde{\d}$ bounds a disk in $H_{\infty}$ by Dehn's lemma,
 and thus $\d$ bounds a disk in $H_n$.
 Thus there is a non-separating sphere in $M_n$,
 contradicting Lemma \ref{irred}.
\end{proof}

\begin{lem} \label{X_n}
 Suppose $n \geq 2k-2$. Then $X_n$ has incompressible boundary.
\end{lem}

\begin{proof}
This is a consequence of Lemma \ref{bdy} and Theorem \ref{handle}.
\end{proof}

\section{Proof that $Y_n$ has incompressible boundary}

In this section, we are under the assumption that
 $\mathcal{E}$ satisfies the $m$-lift condition.
 Recall that $H_n^{\prime} = H_n/D_2^1$.

\begin{lem}\label{nohomotopy}
The curve $\partial E^{n-k+1}$ cannot be isotoped in
 $H_n^{\prime}[\mathcal{E}^{n-k}]$ to intersect $D_1^n$
 fewer than two times.
\end{lem}

\begin{proof}
 We have
\begin{eqnarray*}
\pi_1 H_n^{\prime}[\mathcal{E}^{n-k}]
 &=& < x_1, ..., x_n| w_1, ..., w_{n-k}>\\
 &\cong& <x_1, ..., x_{n-1}| w_1, ..., w_{n-k}> * <x_n>
\end{eqnarray*}
 where $x_j$ is dual to $D_1^j$,
 and the word $w_j$ corresponds to $\partial E^j$.

 The word $w_{n-k+1}$ can be cyclically permuted to have the form
 \begin{equation*}
 w_{n-k+1} = \mathcal{W}_1 x_{n}^{\ell_1}
 \mathcal{W}_2 x_{n}^{\ell_2} ...  \mathcal{W}_t x_{n}^{\ell_t},
\end{equation*}
 where the $\mathcal{W}$'s are freely reduced words
 involving only $x_j$'s with $n-k+1 \leq j \leq n-1$, and each $\mathcal{W}_j$
 represents a nontrivial element in
 the group $< x_1, ..., x_n| w_1, ..., w_{n-k}>$.

 Suppose $\partial E^{n-k+1}$ can be isotoped to be disjoint
 from  $D_1^n$ in $H_n'[\mathcal{E}^{n-k}]$.
 Then, using the relations
 $w_j$, $j \leq n-k$, the word $w_{n-k+1}$ can be re-written entirely in terms
 of $x_j$'s, $j \leq n-1$, and
 this would imply that one of the $\mathcal{W}_j$'s must be trivial
 in $\pi_1  H_n^{\prime}[\mathcal{E}^{n-k}]$.
 However, $\mathcal{W}_j$ is a freely reduced word on
 $x_{n-k+1}, ... , x_{n-1}$,
 which by Theorem \ref{staggered} freely generate a free subgroup,
 for a contradiction.

  Suppose $\partial E^{n-k+1}$ can be isotoped to intersect $D_1^n$
 exactly once.  Then, as in the proof of Lemma \ref{deg2}, $M$ is fibered,
 and so by Corollary \ref{fib}, $\mathcal{E}$ does not satisfy the $m$-lift
 condition, for a contradiction.
\end{proof}

\begin{lem} \label{Y_n}
If  $n \geq m+k-1$,
 then $Y_n$ has incompressible boundary.
\end{lem}

\begin{proof}
Let $n_0=m+k-1$.
We first prove the result in the case where $n=n_0$.\\
\\
\textit{Claim:}
 The manifold $Y_{n_0}$ has incompressible boundary.\\
\\
\textit{Proof of Claim:}  Recall $H_{n_0}^{\prime} = H_{n_0}/D_2^1$.
 By Condition 2 of Definition \ref{m-lift}, the manifold
\begin{eqnarray*}
 H_{n_0}^{\prime} - \partial \mathcal{E}_m^{(m)}
  = H_{n_0}^{\prime} - \partial \mathcal{E}_{m-1}
\end{eqnarray*}
 has compressible boundary, and by
 Condition 1 of Definition \ref{m-lift},
 $H_{n_0}^{\prime} - \partial \mathcal{E}_m$
 has incompressible boundary; therefore, by Theorem \ref{handle},
 the manifold  $(H_{n_0}^{\prime} - \partial \mathcal{E}_{m-1})[E^m]$
 has incompressible boundary.

 Suppose for some $i \in [1, m]$, that
 $H_{n_0}^{\prime}[\bar{\mathcal{E}}_i]- \partial \mathcal{E}_i$
 has incompressible boundary, where
 $\bar{\mathcal{E}}_i = \mathcal{E}_m - \mathcal{E}_i$.
 By Condition 3 of Definition \ref{m-lift}, there is a compression disk for
 $H_{n_0}^{\prime} - \partial \mathcal{E}_{i-1}$
 which is also a compression disk for
 $H_{n_0}^{\prime}[\bar{\mathcal{E}}_i]- \partial \mathcal{E}_{i-1}$.
 Therefore, by Theorem \ref{handle},
 $H_{n_0}^{\prime}[\bar{\mathcal{E}}_{i-1}] - \p \mathcal{E}_{i-1}$
 has incompressible boundary.

 By induction on $i$, it follows that
 $Y_{n_0} = H_{n_0}^{\prime}[\mathcal{E}_m]$
 has incompressible boundary.\\

 Now, suppose $n > n_0$.
 and proceed by induction on $n$.  Suppose $Y_n$ has
 incompressible boundary.
 The manifold $Y_{n+1}$ is obtained from
 $Y_n$ by adding a 1-handle $Z$ to $H_n^{\prime}$, and then attaching
 a 2-handle $E^{n-k+2}$.
 We claim that $(Y_n \cup Z) - \partial E^{n-k+2}$
 has incompressible boundary.

 Suppose otherwise, and let $\Delta$ be a compressing disk.
 Since $Y_n$ has incompressible boundary, the maximal compression body
 for $\partial (Y_n \cup Z)$ has a unique disk system consisting
 only of $D_1^{n+1}$. Therefore if $\Delta$ is
 non-separating, it is isotopic to $D_1^{n+1}$;
 then $\partial E^{n-k+2}$ can be isotoped off of
 $\partial D_1^{n+1}$  in $Y_n \cup Z$, contradicting Lemma \ref{nohomotopy}.

 Suppose $\Delta$ is separating, so it separates off
 a solid torus $V \subset H_{n+1}^{\prime}$,
 with $V \supset D_1^{n+1}$. Since
 $\partial E^{n-k+2} \cap D_1^{n+1} \neq \emptyset$,
 we have $\partial E^{n-k+2} \subset V$.

 So $Y_{n+1}$ contains the punctured lens space $V[E^{n-k+2}]$.
 By Lemma \ref{nohomotopy}, this lens space cannot be $B^3$.
 So $Y_n$ is not irreducible, contradicting Lemma \ref{irred}.
 This completes the proof that $(Y_n \cup Z) - \partial E^{n-k+2}$
 has incompressible boundary.
 Then by Theorem \ref{handle}, $Y_{n+1}$ has incompressible boundary.
\end{proof}

\section{Proof of Theorem \ref{cyclic} and Theorem \ref{surgery}}

\begin{proof} (of Theorem \ref{cyclic})
 Assume $M$ satisfies the hypotheses of Theorem \ref{cyclic}.
 Then by Lemmas \ref{X_n} and \ref{Y_n}, $M_n$ contains
 an incompressible closed surface $S_n$, which is not parallel
 into $\partial M_n$ by Lemma \ref{nonparallel}.
\end{proof}

\begin{proof} (of Theorem \ref{surgery})
Suppose the hypotheses of Theorem \ref{surgery} are satisfied. Let
$b_n$ be the number of boundary components of $M_n$. Then it is easy
to see that $b_n$  is equal to the largest common divisor of $n$ and
$b$ (recall from the proceeding paragraph of Theorem \ref{surgery}
that $b=|\phi(\m)|$ is a finite integer).
 Let each boundary component
of $\p M_n$ have the coordinate basis induced from the basis
$\{\m,\l\}$ of $\p M$. Let $M(b_np/q)$ denote the closed manifold
obtained by Dehn filling each component of $\p M_n$ with slope
$b_np/q$. Then  it is easy to check that $M(np/q)$ is cyclically
covered by $M_n(b_np/q)$. We have shown that $M_n$ contains an
incompressible
 surface $S_n$.
 When  $n \geq width(\lambda)+b$ as well,  there are $b$ successive lifts
 of $\lambda$ contained in $S_n=\partial Y_n$, which implies that
 $S_n$ has  an annular compression to each component of $\p M_n$,
 with slope $0$.  Since $p>1$,
 then by repeatedly applying  Theorem 2.4.3 of \cite{CGLS} $b_n$ times,
 we see that  $S_n$  remains incompressible in $M_n(b_np/q)$.
 Also, by \cite{Sch}, $M_n(b_np/q)$ is irreducible.
 Hence $M(np/q)$ is virtually Haken.
\end{proof}

\section{Brief description of algorithm}

A. Checking that SnapPea presentations are geometric\\

 Given an algebraic word on $x$ and $y$, we need to know
 if it can be represented
 by a simple closed curve in the boundary of a genus 2 handlebody.
 To do this, we attempt to draw a Heegaard diagram.  So we start
 with four disks in the plane, corresponding to $x, x^{-1}, y$ and $y^{-1}$.
 The word $w$ is represented by a collection of edges
 connecting these disks, and we need to find a representation
 which embeds in the plane.

 To program this on the computer, we start placing edges on the graph
 one by one, as indicated by the word $w$.  When placing an edge, the intial
 point is determined from the previous step, but there may be a choice of
 terminal point.  However, it is possible to keep track of the
 choices which are made, and if an impossible situation is arrived at,
 we retreat to the previous choice,
 and change it. In this way, the computer found for every given word,
 a geometric representation.\\
\\
B. Checking that the $m$-lift condition holds\\

 The only non-elementary step in checking the $m$-lift condition
 algorithmically is
 to find, for a given collection of loops in the fundamental
 group of a handlebody $H$, a specific compressing disk for $H$ which
 is disjoint from the loops.  An algorithm for this was given by Whitehead
 (\cite{W}, or see \cite{S2}), which we implemented on GAP.

 Note that Whitehead's algorithm allows one to determine
 the existence of a compressing disk in polynomial time in the
 length of the word; however, to construct the
 disk explicitly requires exponential time.  For our application, we are saved
 this difficulty, since for each compressing disk $\Delta$ we only need
 to compute $[\partial \Delta] \in H_1(\partial H)$.  This allows the
 algorithm to run in polynomial time.

\begin{tabular}{cccc}
Manifold name & width of 2-handle & satisfies $m$-lift for:
                                   & $n$-fold cyclic cover is large for:\\
m006  &        $k=  3    $  & $m=     2 $   &  $n \geq             4$\\
m007  &          $k= 3         $  & $m=    2  $   &  $n \geq           4$\\
m015  &         $k=   4   $  & $m=          2 $   &  $n \geq            6$\\
m017  &          $k= 3     $  & $m=        1  $   &  $n \geq           4$\\
m029  &         $k=  3     $  & $m=        2  $   &  $n \geq           4$\\
m030  &          $k= 3     $  & $m=        2  $   &  $n \geq           4$\\
m032  &          $k= 3   $  & $m=          2  $   &  $n \geq           4$\\
m033  &          $k= 3   $  & $m=          2  $   &  $n \geq           4$\\
m035  &          $k= 3  $  & $m=           1  $   &  $n \geq           4$\\
m037  &          $k=3  $  &  $m=          1   $   &  $n \geq          4$\\
m045  &          $k=3  $  &  $m=           1  $   &  $n \geq           4$\\
m046  &          $k=3  $  &  $m=          1   $   &  $n \geq          4$\\
m053  &          $k= 3 $  &  $m=          1   $   &  $n \geq          4$\\
m054  &          $k= 3 $  &  $m=            2 $   &  $n \geq            4$\\
m058  &          $k= 3 $  &  $m=            1 $   &  $n \geq            4$\\
m059  &          $k= 3 $  &  $m=            1 $   &  $n \geq            4$\\
m061  &          $k= 3 $  &  $m=            2 $   &  $n \geq            4$\\
m062  &          $k= 3 $  &  $m=            2 $   &  $n \geq            4$\\
m066  &          $k= 3 $  &  $m=           1 $   &  $n \geq            4$\\
m067  &          $k= 3 $  &  $m=            1 $   &  $n \geq            4$\\
m073  &          $k= 3 $  &  $m=            1 $   &  $n \geq            4$\\
m074  &          $k= 3 $  &  $m=            2 $   &  $n \geq            4$\\
m076  &          $k= 3 $  &  $m=            1 $   &  $n \geq            4$\\
m077  &          $k= 3$  &  $m=            1  $   &  $n \geq 4$\\
 m079  &          $k= 3$  &  $m=            1  $   &  $n \geq           4$\\
m080  &         $k=  3 $  &  $m=           1  $   &  $n \geq           4$\\
m084  &          $k= 3 $  &  $m=           1  $   &  $n \geq           4$\\
m085  &          $k= 3 $  &  $m=           1  $   &  $n \geq           4$\\
m089  &          $k= 3 $  &  $m=           1  $   &  $n \geq 4$\\
m090  &          $k= 3 $  &  $m=           1  $   &  $n \geq          4$\\
m093  &          $k= 3 $  &  $m=           2  $   &  $n \geq          4$\\
m094  &          $k= 3 $  &  $m=           2  $   &  $n \geq          4$\\
m104  &         $k=  3 $  &  $m=           1  $   &  $n \geq          4$\\
m105  &         $k= 3  $  &  $m=          3   $   &  $n \geq          5$\\
m110  &          $k= 3 $  &  $m=           1  $   &  $n \geq          4$\\
m111  &          $k= 3 $  &  $m=           1  $   &  $n \geq          4$\\

\end{tabular}

\begin{table}
\begin{tabular}{cccc}
Manifold name & width of 2-handle & satisfies $m$-lift for:
                                       & $n$-fold cyclic cover is large for:\\
                                      m137  &         $k=  3 $  &  $m=           2  $   &  $n \geq          4$\\
m139  &          $k= 4 $  &  $m=           3 $   &  $n \geq           6$\\
m148  &         $k=  3 $  &  $m=           1  $   &  $n \geq          4$\\
m149  &         $k=  3 $  &  $m=           1 $   &  $n \geq           4$\\
m202  &         $k=  4 $  &  $m=           2 $   &  $n \geq 6$\\

m203  &         $k=  4 $  &  $m=           1 $   &  $n \geq           6$\\
m208  &         $k=  4 $  &  $m=           1 $   &  $n \geq           6$\\
m249  &         $k=  5 $  &  $m=           4  $   &  $n \geq          8$\\
m259  &         $k=  5 $  &  $m=           3  $   &  $n \geq          8$\\
m260  &         $k=  5 $  &  $m=           3  $   &  $n \geq          8$\\
m261  &          $k= 3 $  &  $m=           1  $   &  $n \geq          4$\\
m262  &         $k=  3 $  &  $m=           1  $   &  $n \geq          4$\\
m285  &         $k=  3 $  &  $m=           1  $   &  $n \geq 4$\\
m286  &         $k=  3 $  &  $m=           1  $   &  $n \geq          4$\\
m287  &         $k= 5 $  &  $m=           5   $   &  $n \geq         9$\\
m288  &         $k= 5 $  &  $m=           3  $   &  $n \geq          8$\\
m292  &         $k= 5 $  &  $m=           3  $   &  $n \geq          8$\\
m319  &         $k= 3 $  &  $m=           1   $   &  $n \geq         4$\\
m320  &         $k= 3 $  &  $m=           1   $   &  $n \geq         4$\\
m328  &         $k= 4 $  &  $m=           1   $   &  $n \geq         6$\\
m329  &         $k= 4 $  &  $m=           2   $   &  $n \geq         6$\\
m340  &         $k= 5 $  &  $m=           1   $   &  $n \geq         8$\\
m357  &         $k= 4 $  &  $m=           2   $   &  $n \geq         6$\\
m366  &         $k= 4 $  &  $m=           1   $   &  $n \geq         6$\\
m388  &         $k= 4 $  &  $m=           1   $   &  $n \geq         6$\\

\end{tabular}
\caption{Data on SnapPea census manifolds }
\end{table}

\vskip1pc
Math Department\\
SUNY Buffalo\\
jdmaster@buffalo.edu\\
\\
Math Department\\
SUNY Buffalo\\
menasco@buffalo.edu\\
\\
Math Department\\
SUNY Buffalo\\
xinzhang@buffalo.edu

\end{document}